\documentclass[11pt]{article}   

\usepackage{geometry}
\pdfoutput=1  

\usepackage{upgreek}
\usepackage{stmaryrd}

\usepackage{flushend}

\usepackage[usenames, dvipsnames]{color}

\usepackage[pdfa,hidelinks]{hyperref}
\usepackage{optidef}
\usepackage{graphics} 
\usepackage{subfigure}
\graphicspath{{figures/}}
\usepackage{amsmath} 
\usepackage{amssymb}  
\usepackage{balance}

\usepackage{verbatim}

\usepackage{textgreek,upgreek,bm}

\makeatletter
\newcommand\gobblepars{%
	\@ifnextchar\par%
	{\expandafter\gobblepars\@gobble}%
	{}}
\makeatother

\def\wham#1{\smallbreak\pagebreak[3]%
	\noindent\textbf{#1}\ \ \gobblepars}

\newcounter{rmnum}
\newenvironment{romannum}{\begin{list}{{\upshape (\roman{rmnum})}}{\usecounter{rmnum}
			\setlength{\leftmargin}{2pt}
			\setlength{\rightmargin}{4pt}
			\setlength{\itemsep}{1pt}
			\setlength{\itemindent}{5pt}
	}}{\end{list}}

\def\bl#1{{\color{blue}#1}}

\newlength{\noteWidth}
\long\def\notes#1{\ifinner
	{\tiny #1}
	\else
	\marginpar{\parbox[t]{\noteWidth}{\raggedright\tiny #1}}
	\fi}

     \def\notes#1{\typeout{See notes!}}
\def\archive#1{}

\def\fl#1{\notes{\bl{FL:  #1}}}

\def\MC{\text{\MC}}

\def\tilalpha{\tilde{\alpha}}

\def\tilpsi{\widetilde{\psi}}

\def\argmin{\mathop{\rm arg\, min}}

\def\Tol{\text{\rm Tol}}

\newcommand{\field}[1]{\mathbb{#1}}

\def\Re{\field{R}}

\def\epsy{\varepsilon}

\def\eqdef{\mathbin{:=}}

\newtheorem{theorem}{Theorem}[section]

\newtheorem{proposition}[theorem]{Proposition}
\newtheorem{lemma}[theorem]{Lemma}

\usepackage{cleveref}

\Crefname{corollary}{Corollary}{Corollaries}
\Crefname{eqnarray}{eq.}{eqs.}
\Crefname{equation}{eq.}{eqs.}

\Crefname{figure}{Fig.}{Figs.}
\Crefname{tabular}{Tab.}{Tabs.}
\Crefname{table}{Tab.}{Tabs.}
\Crefname{lemma}{Lemma}{Lemmas}

\Crefname{theorem}{Thm.}{Thms.}
\Crefname{definition}{Definition}{Definitions}
\Crefname{section}{Section}{Sections}
\Crefname{proposition}{Prop.}{Propositions}
\Crefname{assumption}{Assumption}{Assumptions}
\Crefname{example}{Example}{Examples}

\def\sq{\hbox{\rlap{$\sqcap$}$\sqcup$}}
\def\qed{\ifmmode\sq\else{\unskip\nobreak\hfil
		\penalty50\hskip1em\null\nobreak\hfil\sq
		\parfillskip=0pt\finalhyphendemerits=0\endgraf}\fi}

\def\bfmx{\bfmath{x}}

\def\bfmz{\bfmath{z}}

\def\bfmu{\bfmath{u}}

\def\bfmx{x}
\def\bfmz{z}

\def\bfmu{u}

\def\FRAC#1#2#3{\genfrac{}{}{}{#1}{#2}{#3}}

\def\ddt{{\mathchoice{\FRAC{1}{d}{dt}}%
		{\FRAC{1}{d}{dt}}%
		{\FRAC{3}{d}{dt}}%
		{\FRAC{3}{d}{dt}}}}

\def\ddr{{\mathchoice{\FRAC{1}{d}{dr}}%
		{\FRAC{1}{d}{dr}}%
		{\FRAC{3}{d}{dr}}%
		{\FRAC{3}{d}{dr}}}}

\def\clB{{\cal B}}

\def\clH{{\cal H}}
\def\clJ{{\cal J}}

\def\cX{c_{\text{\tiny\sf X}}}

\def\cdG{c_{{\text{\lower1pt\hbox{d}}}} }

\usepackage{bbold} 

\RenewEnviron{BaseMiniExclam}[7]{%
	\selectConstraintMult{#1}%
	\renewcommand{\localOptimalVariable}{#2}%
	\begin{subequations}
		\ifthenelse{\equal{#7}{b}}{\allowdisplaybreaks}%
		#4
		\begin{alignat}{5}
			\bodyobj{#2}{#3}{#6}{#5}
			\BODY
		\end{alignat}
	\end{subequations}%
	\setStandardMini
}

\def\transpose{{\intercal}}

\def\fee{\upphi}
\def\feeMPC{\fee^{\text{\tiny\sf  MPC}}}
\def\feeMPCQNom{\fee^{\text{\tiny\sf Nom}}}
\def\feeMPCQRob{\fee^{\text{\tiny\sf Rob}}}

\makeatletter
\newcommand{\oset}[3][0ex]{%
	\mathrel{\mathop{#3}\limits^{
			\vbox to#1{\kern-2\ex@
				\hbox{$\scriptstyle#2$}\vss}}}}
\makeatother

\def\xcev{\oset[-.1ex]{\shortleftarrow}{x}}

\def\uH{\underline{H}}

\def\st{\text{\rm s.t.\,}}

\def\ctheta{\check{\theta}}
\usepackage[export]{adjustbox}

\def\MPCt{\tau}
\def\INITt{t_0}

\def\rand{V}
\def\qNomSol{\text{\tiny\sf Nom}}
\def\qRobSol{\text{\tiny\sf Rob}}

\def\tHor{{\hbox{\tiny$\mathcal{T}$}}} 
\def\sHor{{\hbox{\footnotesize$\mathcal{T}$}}} 
\def\Hor{{\mathchoice{\mathcal{T}}{\sHor}{\tHor}{\tHor}}}

\def\Hcev{{\oset[-.2ex]{\shortleftarrow}{H}}}

\def\clHcev{{\oset[-.2ex]{\shortleftarrow}{\mathcal{H}}}}
\def\clCcev{{\oset[-.2ex]{\shortleftarrow}{\mathcal{C}}}}

\def\ccev{{\oset[-.2ex]{\shortleftarrow}{c}}}

\def\uHcev{{\oset[-.2ex]{\shortleftarrow}{\underline{H}}}}
\def\psicev{{\oset[-.2ex]{\shortleftarrow}{\psi}}}
\def\Psicev{{\oset[-.2ex]{\shortleftarrow}{\Psi}}} 

\makeatletter
\def\thanks#1{\protected@xdef\@thanks{\@thanks
		\protect\footnotetext{#1}}}
\makeatother

\title{Model-Free Characterizations of the Hamilton-Jacobi-Bellman Equation 
	and Convex Q-Learning in Continuous Time }

\author{Fan~Lu, Joel~Mathias, Sean~Meyn, and Karanjit~Kalsi
	\thanks{F.~Lu, and S.~Meyn are with the Department of Electrical and Computer Engineering, University of Florida, Gainesville, FL 32611, USA (e-mail: fan.lu@ufl.edu; meyn@ece.ufl.edu).}%
	\thanks{J.~Mathias is with the  School of Electrical, Computer, and Energy Engineering, Arizona State University, Tempe, AZ, 85281, USA (e-mail: Joel.Mathias@asu.edu).}%
	\thanks{K.~Kalsi is with the Pacific Northwest National Laboratory, Richland, WA 99354, USA (e-mail: karanjit.kalsi@pnnl.gov).}%
	\thanks{SM acknowledges support from  NSF
		award  EPCN 1935389,  and from an Inria International Chair, Paris, France.  }    
}

\def\ProofOf#1{\smallbreak\noindent\textit{#1:} \  }

\begin{document}

	\maketitle

	
\begin{abstract}

Convex Q-learning is a recent approach to reinforcement learning, motivated by the possibility of a firmer theory for convergence,  and the possibility of making use of greater a~priori knowledge regarding policy or value function structure.     This paper explores algorithm design in the continuous time domain,   with  finite-horizon optimal control objective. The main contributions are  		
\begin{romannum}
			\item
			Algorithm design is based on a new  \textit{Q-ODE},   which defines the model-free characterization of the Hamilton-Jacobi-Bellman equation.
			
			\item
			The Q-ODE motivates a new formulation of Convex Q-learning that avoids the approximations appearing in prior work.     
The Bellman error used in the algorithm is defined by filtered measurements,  which is beneficial in the presence of measurement noise.

			\item  
			A characterization of boundedness of the constraint region is obtained through a non-trivial extension of recent results from the discrete time setting.
			
			\item   
			The theory is illustrated in application to resource allocation for distributed energy resources, for which the theory is ideally suited.    
		\end{romannum}
		
	\end{abstract}


	


	\section{Introduction}
	
	This paper concerns optimal control of the controlled nonlinear state space model  
	\begin{equation}
		\label{e:dyn}
		\ddt x_t = F(x_t,u_t,t)  \,,\qquad x_0\in \Re^n\,,
	\end{equation}
	in which  the state $x$ and input $u$ evolve on $n$ and $m$ dimensional Euclidean space, respectively.   The goal is to approximate the solution to the finite time-horizon optimal control problem with objective
	\begin{equation}
		\label{e:oc} 
		J(x) =
		\int_{0}^{\Hor} c(x_t, u_t, t)\,dt + J_0(x_\Hor) \,, \ \ x=x_0\,,
	\end{equation}   
	with cost function  $c\colon\Re^n\times\Re^m\times\Re \to\Re_+$,  terminal cost $J_0\colon\Re^n\to\Re_+$,  and fixed $\Hor>0$.    
	
	The minimization over all inputs defines the value function $J^\star(x)$,   and for  $ T_0 \in [0, \Hor)$, the \textit{cost-to-go} is denoted 
	$J^\star(x,T_0)$;  this is the same minimum, but with the integral starting at $T_0$ rather than $0$.   
	Subject to conditions on the model and cost functions, this satisfies the celebrated Hamilton-Jacobi-Bellman (HJB) equation and a characterization of the optimal input as state feedback $		u^\star_t = \fee^\star(x^\star_t,t)$. 
	
	Approximation of a solution is based on concepts from reinforcement learning (RL),  and in particular convex Q-learning.    A starting point is 
	a sample path bound implied by the HJB equation: 
	\begin{equation}
		c (x_t, u_t ,t)  +  \ddt J^\star(x_t, t)  \ge 0\,, \qquad 0\le t\le \Hor
		\label{e:hjbQ}
	\end{equation}
	This holds for any input-state sequence, and is  tight, in the sense that the lower bound is achieved for any $t$ for which $u_t = \fee^\star(x_t,t)$. 
	This inequality could be used in the formulation of a Q-learning algorithm based on linear programming techniques, as in the prior work \cite{mehmey09a,mehmeyneulu21,CSRL,lumehmeyneu22}.   Given a function class $\{ J^\theta  : \theta\in\Re^d\}$ we would use 
	\eqref{e:hjbQ} with $J^\star$ replaced by $J^\theta$ to define a constraint in a nonlinear program, that would be linear if the function class is linearly parameterized.    
	
	The approach proposed in this paper is based on this idea, but with several steps introduced to avoid the use of derivatives of observations.   We arrive at filtering techniques similar to what was introduced in \cite{mehmey09a}, but the finite-horizon setting leads to an exact algorithm that wasn't found in the infinite horizon setting of this prior work.

	\wham{Contributions}
	\begin{romannum}

		\item
		The inequality \eqref{e:hjbQ} is refined to define the \textit{Q-ODE}:    a model-free characterization of the  HJB equation that lends itself to RL algorithm design--see \Cref{t:convexQbdd}.

		\item A new exact formulation of convex Q-learning is obtained, extending \cite{mehmey09a}. The approach is  designed to avoid the numerical challenges that arise in continuous time models.
		
		\item Convex Q-learning is always feasible, but boundedness has been an open topic for research.      Boundedness of the constraint region was characterized  in \cite{lumehmeyneu22} for models in discrete time.  The generalization to the continuous finite time-horizon setting of the present paper is entirely non-trivial--see \Cref{t:AssRel} and \Cref{t:AssConclusions}.
		These results are obtained in the general setting with linear function approximation, so in particular the state space need not be finite.
		
		\item   Results from numerical experiments illustrate several applications of these findings:  (a)  the marriage of MPC and Q-learning is facilitated because constraints can be imposed to ensure convexity of the value function approximation;  (b)   special structure in the application considered,  an extension of \textit{economic dispatch},  justifies a low dimensional function approximation architecture for  convex Q-learning.
	\end{romannum}

\fl{Replace ``\textit{economic dispatch}" with ``allocation for distributed energy resources"?}

	\wham{Related Research}
Convex Q-learning is a recent technique in RL.   While the first paper \cite{mehmey09a} is over one decade old, and did focus on continuous time as in the present paper, this original work laid out theory without attention to algorithms.     The introduction of practical algorithms came only recently in \cite{bascurkraneu21,mehmeyneulu21,lumehmeyneu22}  (see \cite[Ch.~5]{CSRL} for more history,  and \cite{RLandControl13} for a history of RL in continuous time).
	
	Of course, Q-learning has a much longer history.   Watkins' original algorithm \cite{wat89,watday92a} was inspired by older temporal difference learning techniques,   and versions of the temporal difference are also part of convex Q-learning architectures.

The numerical work surveyed in this paper focuses in large part on the marriage of RL and MPC.   The usefulness of an approximate Q-function in MPC was first investigated in  the dissertation \cite{AnuThesis13,kowmaykalmey12} for deterministic control systems,  and contemporaneously in \cite{zhong2013value} for MDPs.    These works are particular approaches to \textit{rollout} for approximate dynamic programming---see \cite{bertsekas2005dynamic} for a survey.    There is also some connection with the older Lyapunov function approach to control design for nonlinear control systems \cite{dixbehdawnag03}.

	Beyond value function approximation, machine learning techniques have been used to approximate unknown or partially known system dynamics in MPC --- see e.g.\  \cite{ostafew2016learning, soloperto2018learning} for use of Gaussian Process regression to learn unmodeled dynamics. Conversely, the use of MPC to shorten the time horizon in entropy regularized RL has been utilized in \cite{bhardwaj2020information, williams2017information}.

	\wham{Organization}   
	\Cref{s:optcontrol} sets the stage, with  a review of MPC, optimality equations, and how these lead to the Q-ODE.
	New Q-learning algorithms are introduced in \Cref{s:algo} based on convex programming,  along with new theory characterizing boundedness of the constraint region.
Application to power systems operations is surveyed in \Cref{s:DD}. 
Conclusions and directions for future research are presented in \Cref{s:conc}.

	\section{HJB Representations}
	\label{s:optcontrol}

	A starting point in the derivation of the HJB equation is Bellman's principle of optimality, which is itself described in terms of the \textit{cost-to-go}:  for each  $ T_0 \in [0, \Hor)$, this is denoted
	\begin{align}
		\begin{split}
			J^\star(x,T_0)  \eqdef  \displaystyle \inf    \Big\{ \int_{T_0}^{\Hor}     c(x_t,& u_t,t) \, dt     + J_0(x_\Hor)   \Big\}
		\end{split}
		\label{e:V0}
	\end{align}
	where the infimum is over continuous $u$ on $[T_0,T]$,   subject to dynamics \eqref{e:dyn},   and   with $x_{T_0} = x$. 	
	
	The principle of optimality is expressed as the family of fixed point equations:   for $\MPCt  \in [0, \Hor)$ and with $x_0=x$,
	\begin{equation*}
		J^\star(x) = \inf_{u_0^\MPCt}    \Big\{  \int_{0}^{\MPCt} c(x_t, u_t, t)\,dt + J^\star(x_\MPCt, \MPCt)  \Big\}
	\end{equation*} 
Dividing each side by $\MPCt$, and letting $\MPCt\downarrow 0$ leads to the  HJB equation:  	
using the shorthand notation  
$J^\star_x = \partial_x J^\star $,  $J^\star_t = \partial_t J^\star $,
\begin{align}
			0 &=  \min_u  Q^\star(x,u,t) 
\label{e:hjb}
\\
 	Q^\star(x,u,t) & \eqdef		  c (x,u,t)  +   J^\star_x(x, t)\cdot F(x,u,t)   +   J_t^\star(x, t)\, .
\nonumber
\end{align}
Letting $\fee^\star(x,t)$ denote the minimizer in \eqref{e:hjb}, an optimal input-state pair is obtained  via state feedback
	\begin{equation*}
		u^\star_t = \fee^\star(x^\star_t,t)
	\end{equation*}

	\wham{Implications to Model Predictive Control} 
 	
We introduce MPC in continuous time only to simplify comparisons to come.     Assumed given is a   ``look-ahead" time horizon $\MPCt$ and ``$\MPCt$-terminal cost'' $c^\bullet$.  For any time $\INITt\ge 0$,  the input $u_{\INITt}$ is obtained through the following steps.  First, the optimization problem is solved: 
	\begin{equation}
		\min\Bigl\{  \int_{t_0}^{t_0+\MPCt} c(x_{\INITt+t}, u_{\INITt+t}, \INITt+t)\,dt + c^\bullet (x_{\INITt+ \MPCt})  \Bigr\}
		\label{e:mpcCC}
	\end{equation} 
	with $x_{\INITt}$ given. The optimizer is a function of time $\{ u^\circ_t : t_0\le t \le t_0 + \MPCt\}$.   The MPC input is defined by $u_{t_0} =  u^\circ_{t_0}$.   
	
	Typically $\MPCt$ will be much smaller than $\Hor$, which in general will lead to performance degradation.   
	However, it follows from the principle of optimality that  the MPC algorithm will minimize the finite horizon objective function \eqref{e:oc} if the $\MPCt$-terminal cost is time varying, with
\[
	c^\bullet (x, \INITt + \MPCt) = J^\star(x, \INITt + \MPCt), \ \ \textit{	for each  $x$, $\MPCt$ and $\INITt$. }
\]
  This ideal is approximated using reinforcement learning techniques in  \cite{AnuThesis13,kowmaykalmey12,zhong2013value}.

	\wham{Q-ODE}

The Q-ODE is  a model-free characterization of the HJB equation \eqref{e:hjb}, inspired by the sample path inequality \eqref{e:hjbQ}.    
	
The function $Q^\star$ that is minimized    in \eqref{e:hjb} is often called the   Q-function,  and easily leads to fixed point equations for reinforcement learning algorithm design for infinite-horizon discounted-cost optimal control  \cite{mehmey09a}.  In this paper RL design is based on approximation of a different function:    fix a scalar  $\sigma>0$ and denote,
\begin{align}  
		H^\star(x&,u, t) \eqdef  - \sigma  J^\star(x, t)  	+ Q^\star(x,u,t)	
\label{e:Hfn}
\end{align}
The use of the letter ``$H$''  recalls the close connection between the Q-function and the Hamiltonian  \cite{mehmey09a}.
The optimal policy is equivalently expressed in terms of this function,
	\begin{equation}
		\fee^\star(x,t) = \argmin_{u} H^\star(x,u,t)
		\label{e:optPolicy}
	\end{equation}

The first step in the Q-ODE construction is the application of the chain rule to obtain, 
\begin{align}
		\begin{split}
			H^\star(x_t,u_t, t) &=   -\sigma   J^\star(x_t, t)  
			\\ 
				& \quad + \bigl[ c (x_t,u_t,t)   + \ddt J^\star(x_t,t)  \bigr]
		\end{split}
\label{e:H_intermediate}
\end{align}
This is valid for \textit{any} input-state trajectory.
The next step is to eliminate $J^\star$ from \eqref{e:H_intermediate}, which requires additional notation.    For any continuous function $H \colon \Re^n\times\Re^m\times \Re \to \Re$, denote $\uH(x,t) = \min_u H(x,u,t)$.   
Application of \eqref{e:hjb} gives 
$
		\uH^\star(x,t) = -\sigma J^\star(x,t)
$,  which on substituting into \eqref{e:H_intermediate} and rearranging terms imples the ODE,
\begin{equation}
\begin{aligned}
			\ddt \uH^\star(x_t, t)  & = \sigma  \uH^\star(x_t, t)  
			\\
			&\quad + \sigma \bigl[ c (x_t,u_t,t)  - H^\star(x_t,u_t, t)   \bigr ]
			\\
\uH^\star(x_\Hor,\Hor) &= -\sigma J^\star(x_\Hor,\Hor) = -\sigma J_0(x_\Hor)	\end{aligned} 		
\label{e:QdvantageCts}
\end{equation}
in which the second equation is treated as a  boundary condition for the first.
This motivates a time-reversal:
For any function $H \colon \Re^n\times\Re^m\times \Re \to \Re$, its time-reversal along an input-state trajectory is denoted  $\Hcev_r \eqdef H\bigl(x_{\Hor - r},u_{\Hor - r}, \Hor-r )$.   When applied to $\uH^\star$, this becomes $	 \uHcev^\star_r = \uH^\star(x_{\Hor - r}, \Hor - r)$.

	
Equation \eqref{e:QdvantageCts} is transformed to the Q-ODE:
	
\wham{Q-ODE}	
With boundary condition $\uHcev^\star_0 = -\sigma J_0(x_\Hor) $,
	\begin{align}
		\label{e:q_ode}
		\ddr \uHcev^\star_r  = -\sigma  \uHcev^\star_r    - \sigma   [ \ccev_r - \Hcev^\star_r  ]    \,,\quad 0\le r\le \Hor\, .
\end{align}

\smallskip

Solutions of \eqref{e:q_ode} involve the  filtered signals,
\begin{subequations}
\begin{align}
\clHcev^\star_r &\eqdef
			\sigma \int_0^r  e^{- \sigma (r-s)}  \Hcev_s^\star   \, ds 
\label{e:smoothedH}
			\\
			\clCcev_r &\eqdef
			\sigma \int_0^r  e^{-\sigma(r-s)}    \ccev_s     \, ds
\label{e:smoothedC}
\end{align}
They are, of course, the solution to the ODEs
\begin{align*}
			\ddr \clHcev^\star_r &=  -\sigma[\clHcev^\star_r   -  \Hcev^\star_r  ]
			\\
			\ddr \clCcev_r &= -\sigma[\clCcev_r  - \ccev_r  ],
		\end{align*}
		with boundary conditions $\clHcev^\star_0 = \clCcev_0 = 0$.
\label{e:smoothed}
\end{subequations}

	The integral representation of the Q-ODE 
	\begin{equation*}
		\uHcev^\star_r = -\sigma e^{ -\sigma r}   J_0(x_\Hor) - \sigma \int_0^r  e^{-\sigma (r-s)}   [ \ccev_s - \Hcev^\star_s  ]      \, ds
	\end{equation*} 
	is thus transformed into the  \textit{algebraic representation}  
	\begin{equation}
		\uHcev^\star_r = -\sigma e^{ -\sigma r}   J_0(x_\Hor) + \clHcev^\star_r - \clCcev_r \,,\quad 0\le r\le \Hor\, .
		\label{e:q_algebraic}
	\end{equation}

	The Q-ODE in the form \eqref{e:q_algebraic} and the proposition that follows will inspire the MPC-Q algorithms surveyed in the next section. 
	The proof of \Cref{t:convexQbdd} is in the Appendix.

	

\begin{proposition}
\label[proposition]{t:convexQbdd}
Suppose that a continuously differentiable solution to the HJB equation exists,  and that an optimal policy is obtained from the minimizer in \eqref{e:hjb}.
		
Suppose  that $H \colon \Re^n\times\Re^m\times \Re \to \Re$ is continuous, and  satisfies the following bound for every $r$ and every input-state trajectory:
		\begin{equation}
			\begin{aligned}
				&   \uHcev_r \ge -\sigma e^{ -\sigma r}   J_0(x_\Hor) + \clHcev_r - \clCcev_r
				\\
				&\textit{with}\ \	\clHcev_r =
				\sigma \int_0^r  e^{- \sigma (r-s)}  \Hcev_s    \,ds .
			\end{aligned}
			\label{e:q_integral_lower}
		\end{equation} 
Then $H(x,u,r)\ge H^\star(x,u,r)$ for all $x,u,r$.
 		\qed
	\end{proposition}

	
	\section{Q-learning Algorithms}
	\label{s:algo}
	
	The algorithms introduced here are based on a  family of approximations $\{ H^\theta : \theta\in\Re^d\}$.      
	For each $\theta$, the \textit{$H^\theta$-greedy policy} is defined by
	\begin{equation}
		\fee^\theta(x,t) = \argmin_u  H^\theta(x,u,t)
		\label{e:Hgreedy}
	\end{equation}
	The ultimate goal of Q-learning is to find the parameter  $\theta^\star$ that leads to the best performance among these policies.   An indirect approach is usually applied, such as the \textit{projected Bellman equation} favored in   much of the academic research.  If we are so fortunate that  $H^{\theta^\star}$ approximately solves \eqref{e:q_algebraic}, then inverse dynamic programming arguments yield  bounds on the performance of the  $\theta^\star$-greedy policy \cite{CSRL}.  
	
	The algorithms described in the following are motivated by \Cref{t:convexQbdd}, which motivates the following definition of the \textit{Bellman error},
	\begin{align}
		\begin{split}
			\clB^\theta_r \eqdef -   \uHcev^\theta_r& - \sigma e^{ -\sigma r}   J_0(x_\Hor) + \clHcev^\theta_r - \clCcev_r .
		\end{split}
		\label{e:q_integralBE}
\end{align} 
in which the filtered signal  $\{ \clHcev^\theta_r : 0\le r\le \Hor\}$ 
is defined as in \eqref{e:smoothedH}.
The inequality \eqref{e:q_integral_lower} using $H=H^\theta$ is equivalently expressed $   \clB^\theta_r\le 0 $ for each $r\in[0,\Hor]$.

	\wham{Projected Bellman Error}  
	
The algorithm described here is inspired by the  DQN algorithm as described in  \cite{CSRL}.  
	\begin{subequations}

		\fl{I changed 'frac' to 'tfrac' since 'frac' is too big and the bracket is too small.}
		For each $n$, given the current estimate $\theta_n$,  the parameter update is obtained as the solution to the nonlinear program, 
		\begin{align}
			&
			\theta_{n+1} = \argmin_\theta \bigl\{ \|     \clB^{\theta\mid\theta_n}  \|_{L_2}^2  
			+\tfrac{1}{\alpha_{n+1}} \| \theta - \theta_n\|^2\bigr\}
			\label{e:theta_recur}
			\\[.5em]
			&
			\begin{aligned}
				\clB_r^{\theta\mid\theta_n} \eqdef -   \uHcev_r^{\theta_n} &  - \sigma e^{ -\sigma r}   J_0(x_\Hor) 
			+ \clHcev^\theta_r - \clCcev_r
			\end{aligned}
			\label{e:DQN_Bell1}
		\end{align}
		in which the non-negative  sequence $\{\alpha_n : n\ge 1\}$ is analogous to the usual step-size sequence in RL.   The term \eqref{e:DQN_Bell1}  is defined as in \eqref{e:q_integralBE}, with the first appearance of $\theta$ frozen.    
		The $L_2$ norm in \eqref{e:theta_recur} is the standard, $ \|     \clB^{\theta\mid\theta_n}  \|_{L_2}^2  \eqdef \int_0^\Hor  [    \clB^{\theta\mid\theta_n}  _r]^2  \,  dr$.
		
		\label{e:ourDQN}
	\end{subequations}
	
	The algorithm is simplified significantly when the parameterization is linear:  
	\begin{equation}
		\label{e:linearpara}
		H^\theta(x,u,r) =   \theta^\intercal \psi(x,u,r),
	\end{equation}
	where the $d$-dimensional basis $\psi$ might be chosen based on known structure of the control problem.  In this case, we write	\begin{equation}
		\label{e:slPsi}
		\Psicev_r \eqdef \sigma  \int_0^r  e^{- \sigma (r-s)} \psicev_s   \, ds \,,\quad 0\le r\le\Hor\,,
	\end{equation}
	with $\psicev_s \eqdef \psi(x_{\Hor - s}, u_{\Hor -s}, \Hor - s)$. This gives  $\clHcev^\theta_r = \theta^\intercal \Psicev_r$, and   \eqref{e:DQN_Bell1} becomes   
	\begin{align}
		\label{e:Bell_algebraic}
		\begin{split}
			\clB^{\theta\mid\theta_n}_r =  -   \uHcev_r^{\theta_n}- \sigma e^{ -\sigma r}   J_0(x_\Hor) 
			  + \theta^\transpose \Psicev_r  -  \clCcev_r 
		\end{split}
	\end{align}

	Substituting \eqref{e:Bell_algebraic} in \eqref{e:theta_recur} and taking the gradient with respect to $\theta$ leads to the fixed point equation that is solved to obtain $\theta_{n+1} $,
	\begin{equation}
		0=	\langle \clB^{\theta_{n+1} \mid\theta_n}  ,  \Psicev \rangle   +  \tfrac{1}{\alpha_{n+1}} [ \theta_{n+1} - \theta_n ] 
		\label{e:DQN}
	\end{equation}
in which  the first term depends linearly on $\theta_{n+1}$:
\begin{align*}
		\langle \clB^{\theta_{n+1} \mid\theta_n}  ,  \Psicev \rangle \eqdef \int_0^\Hor \Psicev_r \clB^{\theta_{n+1} \mid\theta_n} _r\,dr \, .
\end{align*}
	  	
	If the resulting sequence of estimates $\{\theta_n\} $ is bounded,   it follows that $\| \theta_{n+1} -\theta_n\| = O(\alpha_{n+1})$,  which justifies the following \textit{approximation}:  
	\begin{equation}
		\theta_{n+1} = \theta_n - \alpha_{n+1} \langle \clB^{\theta_n \mid\theta_n}  ,   \Psicev \rangle 
		\label{e:batchQ}
	\end{equation}
	This is a variation of Watkins'  algorithm, with two significant changes:  The basis is not tabular (so   the function class does not span all possible functions of state and action),    and the temporal difference in  \cite{wat89,watday92a}  is replaced by the special version of the Bellman error introduced in this paper.
	
	There is currently no theory to predict the success of DQN \eqref{e:DQN} or the recursion \eqref{e:batchQ}.    Stability of Q-learning is largely an open topic outside of very special cases  (see discussion in  \cite[Section 3.3.2]{sze10},    \cite[Section 11.2]{sutbar18} and \cite[Chs.~5,9]{CSRL}). 
	
\wham{Convex Q-Learning}
 \Cref{t:convexQbdd} is motivation for the following  ``ideal'' algorithm:
 Choose a probability measure $\mu$ on $\Re^n\times\Re^m\times[0,\Hor]$, and solve the nonlinear program,
	\begin{subequations}
		\begin{align} 
			\theta^\star = \argmin_\theta
			\ \  &   \langle \mu, H^\theta \rangle     
			\\
			\st  \ \   
			& 
			{ \clB^\theta_r }    \le 0, \qquad r \in [0, \Hor]
			\label{e:LPQ1constraint_cts}
		\end{align}
		
		\fl{note also here $J_0$
		\\
		I don't understand what you mean}

		\label{e:LPQ1reg_cts}%
	\end{subequations}%
	
	In practice, the infinite number of constraints in
	\eqref{e:LPQ1constraint_cts}
	must be relaxed.  In this paper, we replace the constraints in \eqref{e:LPQ1constraint_cts}
	by the single constraint,
	\begin{equation}
		\frac{1}{\Hor}
		\int_0^{\Hor} \bigl[ \clB^\theta_r \bigr]_+\, dr  \le \Tol
		\label{e:relCons}
	\end{equation}
	where $\Tol>0$ is a small constant,  and $[s]_+=\max(0,s)$.
	
	This is a convex program under mild conditions.    The convexity assumption in \Cref{t:ConvexQ} will hold when the function class is linear (i.e., defined with respect to a basis via  \eqref{e:linearpara}).  
	\begin{proposition}
		\label[proposition]{t:ConvexQ}
		Suppose that the Bellman error \eqref{e:q_integralBE} is a convex function of $\theta$ for each $r$.   Then the constraint regions 
		\eqref{e:LPQ1constraint_cts} and \eqref{e:relCons} are each convex subsets of $\Re^d$.   	 
		\qed
	\end{proposition}

	\subsection{Exploration and Constraint Geometry}
	
	The following assumptions are imposed throughout the remainder of the paper:

	\wham{Assumption A1:}
	The function class is linear,  $\{ H^\theta = \theta^\transpose \psi : \theta\in\Re^d\}$.
	The basis $\psi\colon\Re^n\times\Re^m\times \Re_+\to\Re^d$ and the cost function  $c\colon\Re^n\times\Re^m\times \Re_+\to\Re_+$ are continuously differentiable ($C^1$). 
	
	Moreover, for each $\theta\in\Re^d$,  the minimum in \eqref{e:Hgreedy} defines a continuous feedback law $\fee^\theta \colon\Re^n  \times \Re_+ \to\Re^m$.    And,  with  $u_t = \fee^\theta(x_t,t)$ for $0\le t\le\Hor$ there is a solution to the state equation \eqref{e:dyn}.  
	\qed

The constraint set associated with 	  \eqref{e:relCons} is denoted
	\begin{equation}
		\Uptheta =\Bigl\{ \theta\in\Re^d :   
		\frac{1}{\Hor}\int_{0}^{\Hor} \bigl[\clB^\theta_r\bigr]_+\, dr \le \Tol  \Bigr\}
		\label{e:thetaConstraints}
\end{equation}
It is always non-empty since it contains the origin,  but boundedness of $\Uptheta$ has been an open topic for research.  

Necessary and sufficient conditions for boundedness will be obtained based on  algebraic conditions on the basis along input-output sample paths obtained for training.  To ease analysis and save space, we adopt the  notation,  
\[
\psi_t \eqdef \psi(x_t,u_t,t)\,, \qquad
	\psicev_r \eqdef \psi_{\Hor - r} \,,  \quad 0\le t, r\le \Hor\,.
\]
The covariance matrix is denoted
	\begin{equation}
		\label{e:covama}
		\Sigma \eqdef \frac{1}{\Hor}\int_0^\Hor \tilpsi_s\tilpsi_s^\transpose \, ds\,,
		\quad \textit{with}  \ \ \tilpsi_s \eqdef \psi_s  - \frac{1}{\Hor} \int_0^\Hor \psi_t \,  dt
	\end{equation}
 	
The conditions that follow are the focus of analysis in the remainder of this section. The third is a standard assumption intended to capture 
``sufficient exploration'' in temporal difference learning \cite{tsivan97,CSRL}.     In the context of this paper, it is Condition~E1 that is most valuable:
\Cref{t:AssConclusions} tells us  that $\Uptheta$ is bounded under this condition, and hence 
what should be considered ``good exploration''.


	\begin{romannum}
		\item[\textbf{Condition E1:}]    The set $\{ \psi_t : 0 \le t \le \Hor \}$ is not restricted to any half space in $\Re^d$.
		
		\item[\textbf{Condition E2:}]  
		The only vector $v\in\Re^d$  satisfying $\uHcev^v_r \ge \clHcev^v_r$ for all $0 \le r \le \Hor$ is $v=0$.  
		
		\item[\textbf{Condition E3:}]  $\Sigma>0$, with $\Sigma$ defined in \eqref{e:covama}.
	\end{romannum}

	\begin{proposition}
		\label[proposition]{t:AssRel}
		If Condition E1 holds then Conditions E2 and E3 follow.  
	\end{proposition}
We postpone the details of the proof to the Appendix.   
The relationship between E1 and E3 is straightforward, since the latter is equivalent to the statement that $\{ \psi_t : 0 \le t \le \Hor \}$ is not restricted to any \textit{subspace} in $\Re^d$.     As part of the proof that E1 implies E3,  it is shown that if $\{\psi_t\}$ is not restricted to any half space in $\Re$, so is the difference $\{\psicev_r - \Psicev_r\}$.


\Cref{t:AssRel} combined with the following establishes that boundedness of 	$\Uptheta$ is equivalent to Condition~E2.  The proof is postponed to the Appendix.    

\begin{proposition}
		\label[proposition]{t:AssConclusions}
		If Condition~E1 holds then $\Uptheta$ is bounded.     
		Conversely, if $\Uptheta$ is bounded then Condition  E2  holds.
	\end{proposition}

\section{Optimal  Dispatch of Energy Resources}
\label{s:DD}

This section serves to illustrate the   marriage of convex Q-learning with MPC,  and show that it may provide efficient solutions to complex control problems found in power systems applications.

We consider the optimal allocation of distributed energy resources (DERs) in a dynamic setting.  
The goal is to schedule  generation and other ``balancing assets'' to meet supply-demand constraints while minimizing cost, similar to economic dispatch.
It was discovered recently that a form of state space collapse can be expected \cite{matmoymeywar19}.    

It is assumed that the balancing assets are derived from   flexible loads (such as water heaters or water pumping) along side batteries.  We will use the term  \textit{virtual energy storage} (VES) for both real and virtual batteries.


\subsection{Dispatch model}  
It is assumed that there are  $M\ge 2$ classes of VES,  in addition to generation (which is modeled as a single resource, i.e., the aggregation of all the traditional generators in the balancing area). The goal is to optimally allocate these resources to balance the net load $\ell$ over the time horizon $[0, \Hor]$.


Following \cite{haosanpoovin15, cammatkiebusmey18,matmoymeywar19},  the \textit{state of charge}  (SoC) for the $i$th VES class is assumed to evolve according to the linear dynamics 
\begin{equation}
	\ddt x_t^i   = -\alpha_i x_t^i    - z_t^i    \qquad 1\le i\le M\,,
	\label{e:SoC_ODE}
\end{equation}
in which     $-z_t^i$ is power deviation at time $t$, and $\alpha_i$ is a non-negative leakage parameter.    For a TCL, the SoC $x_t^i $ is an affine function of internal temperature, and $\alpha_i$  corresponds to the thermal time constant.

\wham{Formulation of a cost function}
 
A cost function is designed based on three goals:  maintain the SoC within bounds, and penalize peaks and ramps in generation. 
 To impose a cost on ramping it is necessary to augment the state description, introducing   
\begin{equation}
	u_t^i \eqdef \ddt z_t^i  \, .
	\label{e:ramp_nu}
\end{equation}
We view \eqref{e:SoC_ODE} and \eqref{e:ramp_nu} as a linear dynamical system with the augmented state $x^a \eqdef (\bfmx, \bfmz)$,
 and control input $\bfmu$.

\begin{subequations}	

The dispatch problem  is  formulated as a
	finite-horizon optimal control problem:
	\begin{align}
		\!\!\!\!\!
		\min  \ \ &  
		\int_0^{\Hor}  c(x_t, z_t, u_t, t)   \, dt  +	J_0(x^a_\Hor)  
		\label{e:qp19}
		\\
		\st  \ \ & 
		\ell_t=g_t+z^\sigma _t
		\label{e:balancecons}
		\\
		&\ddt g_t=\gamma_t
		\label{e:genrampcons}
		\\
		&	\ddt {x}_t^i=- \alpha_i x_t^i - z_t^i
		\label{e:soccons}
		\\
		&	\ddt z_t^i= u_t^i, \ \ 1\le i \le M  \,,\  0\le t\le \Hor
		\label{e:loadrampcons}
	\end{align} 
	with $x^a = (x_0, z_0) \in \Re ^M \text{ given}$,   $z^\sigma_t = \sum z_t^i$.
\end{subequations}	

The constraint \eqref{e:balancecons} ensures that the supply from generation, batteries, and VES matches net load. The dynamics of generation ramping are given by \eqref{e:genrampcons}.

The terminal cost $J_0$ in \eqref{e:qp19} was chosen to be quadratic function, of the form 
$J_0(x,z) = x^\transpose D x + k_\ell(z^\sigma - \ell_\Hor )^2$ with $k_\ell>0$ and $D>0$ diagonal  ($M\times M$).
The cost function $c$  was taken as the sum of three components, reflecting the three goals:  
\begin{align*}
		c(x_t, z_t, u_t, t)  =   \cX(x_t) +  \kappa \bigl[   u^\sigma_t - \ddt\ell_t \bigr]^2 +  \kappa_\ell [ z^\sigma_t- \ell_t]^2 
\end{align*}
with $u^\sigma_t = \sum u_t^i$,   and $ \kappa, \kappa_\ell$ positive constants. 
A soft constraint on capacity is imposed via
\fl{change subscript $i$ to superscript.}
\begin{equation}
	\cX(x) =  \sum_{i=1}^{M}   c^i(x^i)  \,, \quad x\in\Re^M\,,
	\label{e:cX}
\end{equation}
where each $c^i\colon\Re\to\Re_+$ is smooth and strongly convex.   In the numerical results they were chosen convex, of the same form as in a portion of the numerical results from \cite{matmoymeywar19}.



This optimal control problem falls in the category of singular optimal control because the cost is not coercive in $u$   (there  is a cost on the  sum  $u^\sigma_t$, and not on the individual terms $u_t^i$) \cite{fra79, hausil83}.

%

A major conclusion of  \cite{matmoymeywar19} is that the cost to go for any time $T_0$  can be expressed as a convex function of $x^{\sigma,a} \eqdef (x^\sigma, z^\sigma)$  with $x^\sigma = \sum x^i$ (this is a consequence of the state space collapse referred to earlier).  
However, in the prior work  it is assumed that $J_0\equiv 0$.     The conclusions will change when $J_0$ is coercive, as assumed in the numerical results that follow.  However, state space collapse provides ample motivation for the choice of function class in Q-learning.

 \wham{Function approximation architecture}
 
The function approximations considered so far are model based, in which we begin with an affine function class for approximation of the value function: 
\begin{equation}
J^\theta(x^a,t) =   J_0(x^a) + \theta^\transpose \psi(x^{\sigma,a},t)\,,\quad \theta\in\Re^d \,, 
\label{e:Jfn}
\end{equation}
in which $\psi\colon \Re^2\times\Re_+\to\Re^d$.     The representation \eqref{e:Hfn} then motivates the function class, with candidate approximations 
\fl{$Q^\theta$ is a function of $x^a$. The second equation here is exceeding the half page.}
\begin{align*}  
H^\theta(x^a,u,  t)  & \eqdef   - \sigma  J^\theta(x^a, t)     +   Q^\theta(x^a,u,t)
		\\
Q^\theta (x^a,u,t)  &\eqdef		 c (x^a,u,t) 
		 \\
		 &\qquad  +   J^\theta_x (x^a, t)\!\cdot\! F(x^a,u,t)   +J^\theta_t(x^a, t) 
\end{align*}
The vector field $F$ is not difficult to estimate in this particular example .
The impact of model uncertainty is investigated in the numerical results that follow.  

\notes{error in $J_0$ to be fixed in the new year (it must be coercive, and based on Sept. 21 discussions it seems it wasn't.)
			\\
			\fl{ok.}
}

 To match the ideal $J^\theta(x^a,\Hor ) =   J_0(x^a) $,  the basis was designed to ensure $ \psi(x^{\sigma,a},\Hor) =0$ for each $x^{\sigma,a}\in\Re^2$. It is convenient to take a typical basis function of the form   
\begin{equation}
	\psi_{i, j}(x^{\sigma,a}, t) = q_i (x^{\sigma,a}) p_j  (t) 
	\label{e:psiED}
\end{equation}
in which     $q_i \in  \{(x^{\sigma})^2, x^\sigma, (z^{\sigma})^2, z^\sigma, 2x^\sigma z^\sigma , 1 \}$ for $1\le i\le 6$.    The functions $\{ p_j \}$ were taken to be a mixture of Fourier basis elements and polynomials.   Through trial and error we arrived at three possibilities:   we took $p_1 (t) =t^2$,  and for $j \ge 2$ the function $p_j$ was an element of the set 
\[
  \{   
 	1-	\cos(\omega_i t )   :   1\le i \le n_f  \}   
 \]
with $n_f =  30 $ in all experiments.  Thus, $d =5\times 31 = 155$.

The basis was chosen so that the functions of time are non-negative.   Writing $\theta\in\Re^d$ in compatible form so that $\theta^\transpose \psi = \sum_{i,j} \theta_{i, j}	\psi_{i, j}$,   the constraint $\theta_{i,j}\ge 0$ was imposed in implementations of  convex Q-learning for any $i,j$  for which $	\psi_{i, j}(x^{\sigma,a}, t)    = (x^{\sigma}) ^2 p_j (t) $ or $(z^{\sigma}) ^2 p_j (t) $.    It was found that this helped to ensure that the solution $\theta^*$ would result in a cost to go approximation $J^{\theta^*}(x^a, t)$ that is convex in its first variable for each $t$.



\subsection{Simulations}

The system parameters for VES and net load $\ell$ were taken from \cite{matmoymeywar19}. 
The optimal dispatch problem  \eqref{e:qp19} was considered with  $M=5$ VES classes: ACs, residential WHs (fwh), commercial WHs  (swh), refrigerators (rfg), and pool pumps (pp).     
The time horizon $\Hor$ was set to 24 hours,   and $\sigma = 5\times 10^{-4}$ in the convex Q learning algorithm and the definition of $\clH^\theta$.


\wham{Training architectures}
Two strategies were employed to construct the convex program \eqref{e:LPQ1reg_cts}.    In each case data was collected from 44 independent runs, differentiated as explained in the following.

\textbf{1. Q nominal training}:    initial conditions $x^a$ sampled uniformly from $\Re^{2M}$ at random, trajectories are generated from the nominal model with parameter $\alpha\in\Re^M_+$. These trajectories were then used to solve the convex program \eqref{e:LPQ1reg_cts}. The resulting solution of \eqref{e:LPQ1reg_cts} gives a value function approximation denoted $J^{\qNomSol}$  (of the form 
\eqref{e:Jfn} for the parameter estimate $\theta^{\qNomSol}$).

\fl{change subscript to superscript.}
\textbf{2. Q robust training}:  In addition to sampling    initial conditions $x^a$,   in each batch  the model is perturbed via  
\begin{align}
	\tilalpha^i(\epsy, \rand) = \alpha^i \times \rand^i
	\label{e:perturbedM}
\end{align}
where $\{\rand^i\}$ were selected i.i.d.\ and  sampled independently of $x^a$  from $[1 - \epsy, 1+\epsy]$,
 with $\epsy$ ranging from 0 to 1.     
We then generate trajectories with initial condition $x^a$   using the perturbed model $\tilalpha$. 
These trajectories are used to solve \eqref{e:LPQ1reg_cts},  which defines a value function approximation denoted $J^{\qRobSol}$
 (of the form 
\eqref{e:Jfn} for the parameter estimate $\theta^{\qRobSol}$).

\wham{Performance evaluation}  

To evaluate the outcome of convex Q-learning training required additional experiments.   For testing performance in MPC, we note that the policy defined by   \eqref{e:mpcCC} can be defined as (time varying) state feedback.  The policy
$\fee=\feeMPCQNom$ and $\fee=\feeMPCQRob$ were obtained on replacing $c^\bullet$ in \eqref{e:mpcCC} with $J^{\qNomSol}$ and $J^{\qRobSol}$, and $\fee = \feeMPC$ based on MPC with zero penalty term $c^\bullet$.  
The policies based on Q-learning will be called MPC-Q.

For any feedback policy $u_t = \fee(x^a_t, t)$, the associated total cost is denoted
\begin{equation}
J^\fee(x^a) = \int_{0}^\Hor c(x^a_t , u_t, t)\, dt \,, \qquad x^a = (x(0),z(0))
\label{e:JfeeMPC}
\end{equation}
In the numerical results summarized below we compared this with the optimal $J^\star(x^a)$ from specific initial conditions, and also the cost to go.      This section concludes with experiments illustrating the impact of model error,  for which the performance metric was the average over 
independent trials,   with both the initial condition and model perturbed in each trial.

\begin{figure}[ht]
	\centering
	\includegraphics[width=0.6\hsize]{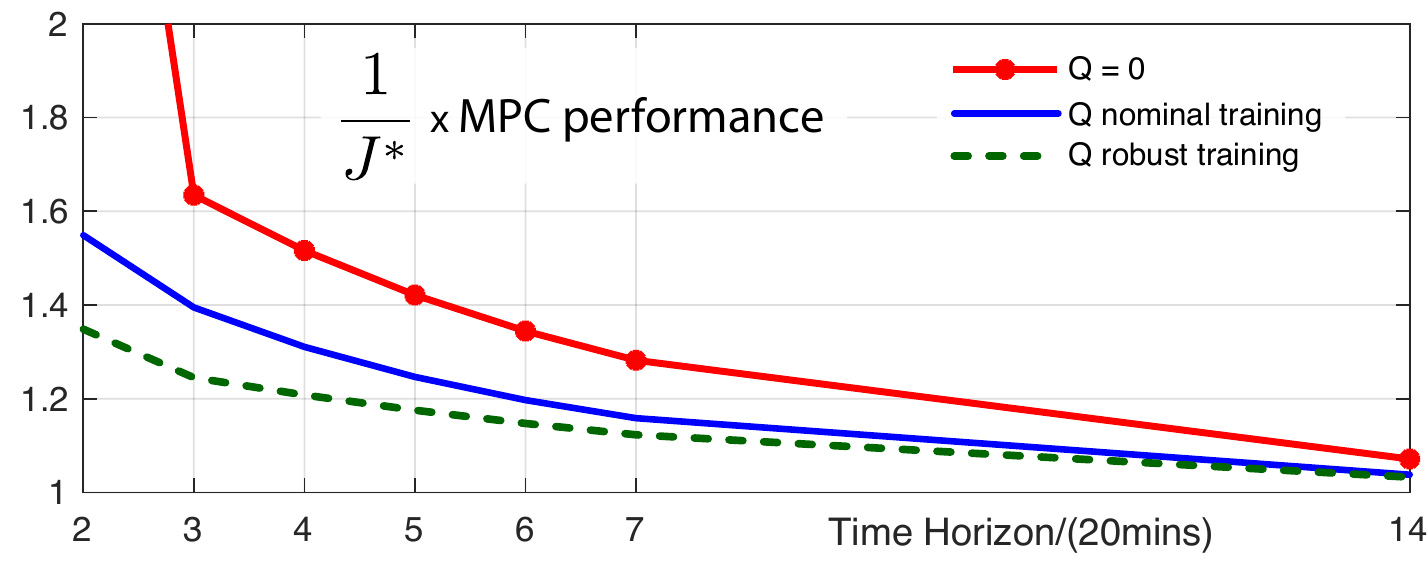}
	\caption{Performance of MPC-Q on the nominal model.
 	}
	\label{fig:totalCost}
\end{figure}

\notes{ I have no idea what this means:
\\
    We performed experiments as follows:  an initial condition $x^a$ was chosen uniformly at random from $\Re^{2M}$,  and the optimal total cost
$J^{\theta^\star}(x^a)$  was compared with the total cost obtained from both MPC-Q and MPC for a range of look-ahead horizons.
\\
It sounds like you simply looked at one initial condition
\\
\fl{Yes, we only look at one initial condition that was different from training. What I mean by comparison is that I was comparing $J^\star$ with $J^\fee$, where $\fee$ was chosen to be $\feeMPC$, $\feeMPCQNom$ and $\feeMPCQRob$ in \Cref{fig:totalCost}}.
\\
We have to discuss this in October.
}

\wham{Experimental results 1:  testing on nominal model.}

The first experiment was designed to investigate the loss in performance introduced from perturbations of the model during training.


\begin{figure}[ht]
	\centering
	\includegraphics[width=0.6\hsize]{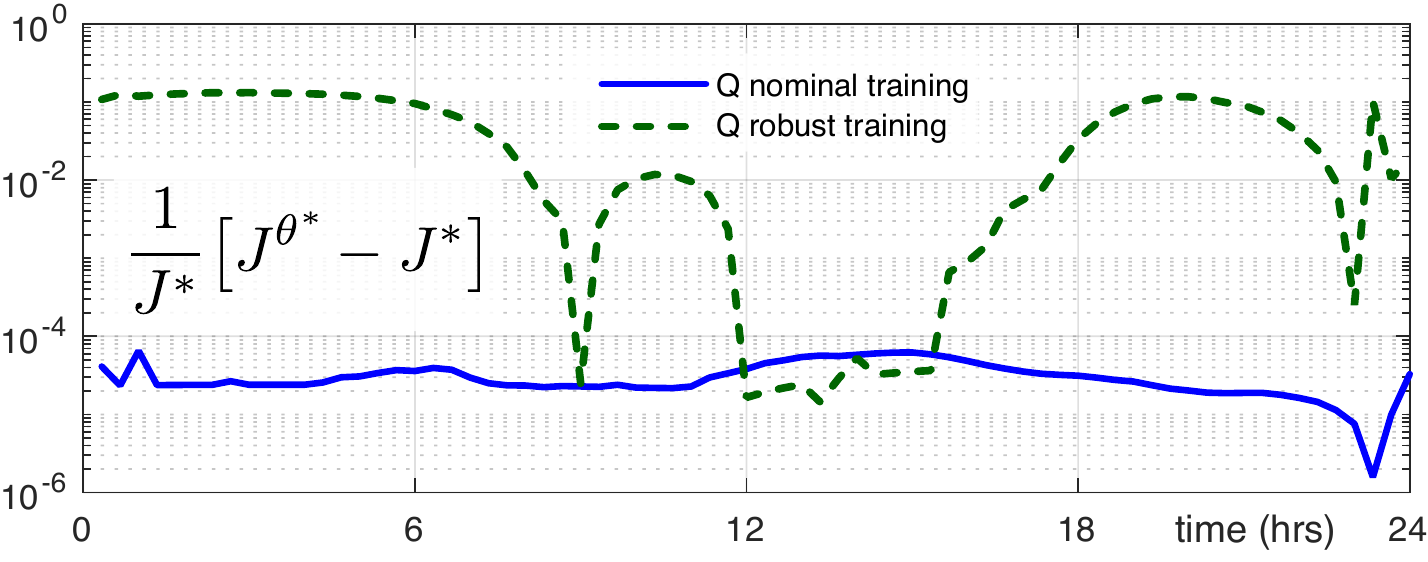}
	\caption{Normalized error with convex Q-learning,  evaluated along an optimal trajectory.   }
	\label{fig:cost2go}
\end{figure}

The normalized error between  the approximation $J^{\theta^*}$ and the optimal cost to go  $J^\star  $ was obtained for the two training approaches with $J^{\theta^*}$ indicating either $J^{\qNomSol}$ (nominal training) or $J^{\qRobSol}$ (robust training).
\notes{I refuse to say $x^a$ was picked at random!
\\
\fl{You're right. For testing on nominal model, we only use one initial condition that is different from training.}
}

In these experiments the initial condition was fixed at a typical value, and   the true optimal solution $\{ x^\star_t, z^\star_t, u^\star_t : 0\le t\le \Hor \}$ was obtained.     For each $t$,   the cost-to-go $J^\star (x^\star_t, z^\star_t, t ) $, was compared with $J^{\qNomSol}(x^\star_t, z^\star_t,  t )$ and $J^{\qRobSol}(x^\star_t, z^\star_t, t )$.  It is seen in \Cref{fig:cost2go}  that  the performance gap using $J^{\qNomSol}$ is less than 0.008\%\ throughout the run. Though the gap using $J^{\qRobSol}$ is larger, we will see its more robust to model perturbations.


\notes{Old caption:  Normalized $J^\fee$ as a function of look-ahead horizon $\MPCt$.
\\
Does this reflect the data ?
\\
\fl{what do you mean by the data?}
\\
Let me think it over.
}

The next results show the performance of MPC as described above \eqref{e:JfeeMPC}. 
\Cref{fig:totalCost} shows data from one typical experiment, performed on the nominal model.   It is surprising to see that
the policy   
$\feeMPCQRob$ gave the smallest error from  $J^\star$ for each look-ahead horizon considered (as small as 40 minutes).    Performance for $\feeMPC$ (with $c^\bullet \equiv 0$) was far worse.  

\begin{figure}[ht]
	\centering
	\includegraphics[width=0.6\hsize]{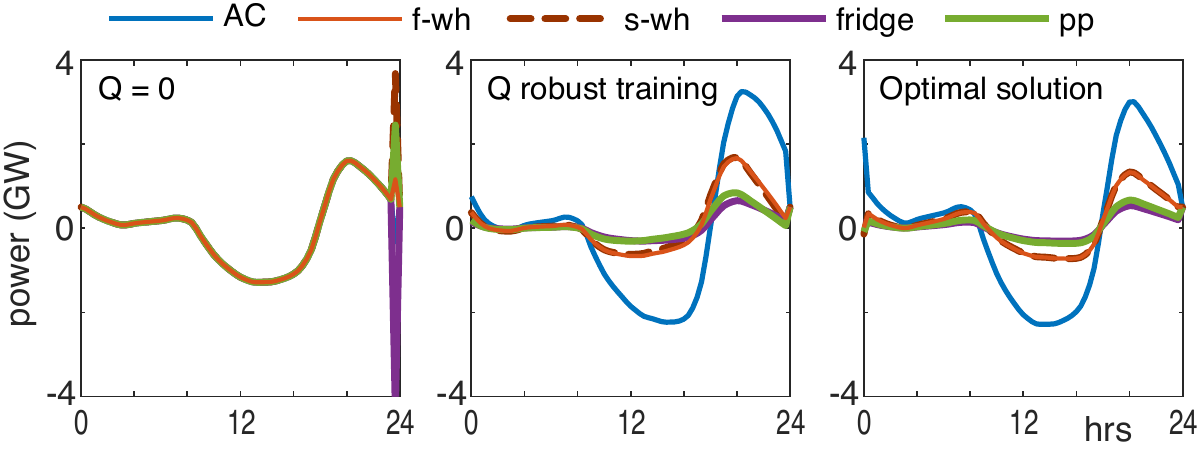}
	\caption{Trajectories of power deviation from each load class.}
	\label{fig:powerTraj}
\end{figure}

\Cref{fig:powerTraj}  shows the power trajectories obtained using MPC-Q with Q robust training mirrors the optimal solution with look-ahead horizon $\MPCt = 40$mins, and how MPC dramatically fails (we omit plots for MPC-Q with Q nominal training since it has similar performance).

\notes{ what does this mean?
\\
 In this experiment the initial condition was chosen uniformly at random from $\Re^M$,
 \\
 Just one initial condition?  Then, what is the point?
	\\
	\fl{I should have plotted \Cref{fig:totalCost} and \Cref{fig:cost2go} using multiple ICs uniformly sampled from the state space and do the average.} 
	\\
	More to discuss in October!
}

\wham{Experimental results 2:  testing on perturbed models}

The impact of model uncertainty is investigated next.

To test a given policy $\fee$ we conducted $N_p$ independent trials for a range of $\epsy\ge 0$,   and averaged the resulting total cost obtained in each trial to obtain 
\begin{equation}
\widehat{J}^\fee_\epsy  = \frac{1}{N_p} \sum_{k=1}^{N_p} J^\fee(x_k^a)  
\label{e:avgPerformance}
\end{equation}
For each $k$ the initial condition was chosen randomly, as well as the perturbation of the model defined by $\tilalpha_k$ via \eqref{e:perturbedM}
for   $1\le k\le N_p$,  with $N_p=50$.   
\fl{change superscript $k$ to subscript.}

\begin{figure}[ht]
	\centering
	\includegraphics[width=0.6\hsize]{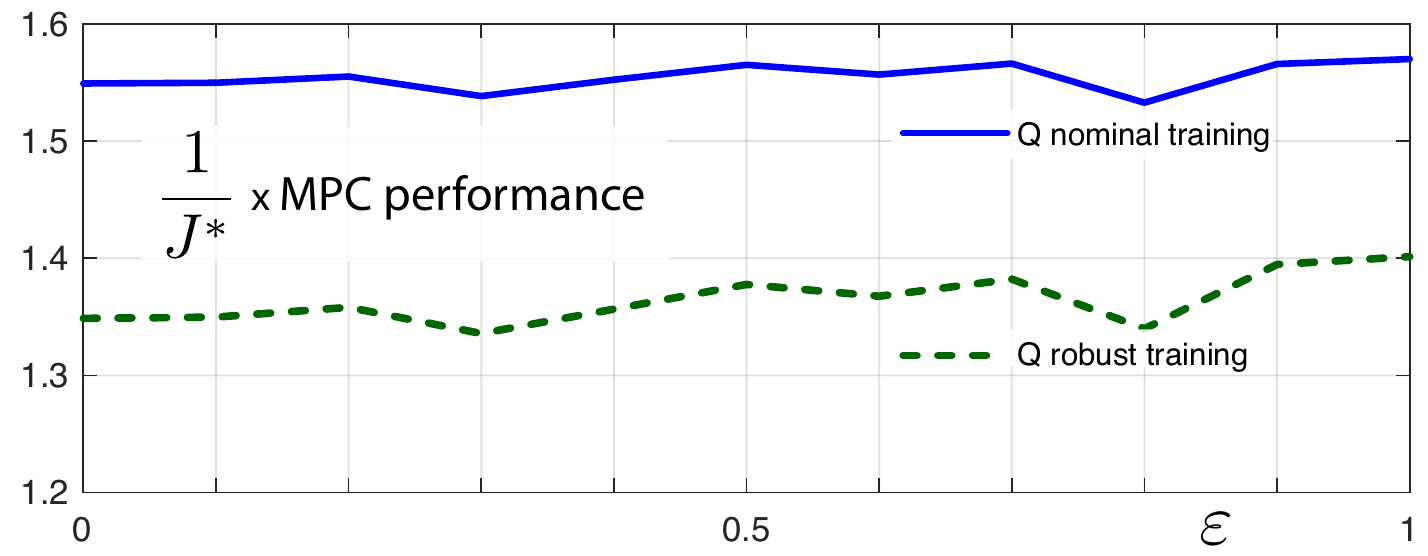}
	\caption{Robustness of MPC-Q: normalized averaged total cost in MPC-Q as a function of $\epsy$ with look-ahead horizon 40mins.}
	\label{fig:errTotalCostFuncEpsy}
\end{figure}

\notes{
	Fan, you had written ``predict'' horizon in the caption.   We are using ``look-ahead'' elsewhere.   This is what you mean, right?
	\\
	\fl{Yes.}
}

 \Cref{fig:errTotalCostFuncEpsy} shows that $\widehat{J}^\fee_\epsy$ is nearly independent of $\epsy$ for either policy $\fee = \feeMPCQNom$ or $\fee = \feeMPCQRob$,  with the latter giving better performance for each value of $\epsy$ tested.   The results for MPC without penalty term is not shown since the ratio was always greater than 3.
 
In all experiments conducted to date,  we find that  Q robust training increases the robustness of MPC-Q.   
Recall that results in \Cref{fig:totalCost} indicate that the robust training can improve closed loop performance even for the nominal model.

\section{Conclusions}
\label{s:conc}
	 
	The Q-ODE  for finite-horizon optimal control is a new model-free characterization of the HJB equation that lends itself to the formulation of reinforcement learning algorithms.  Two were highlighted in this paper: the convex Q-learning algorithm that has been the focus, and a variation of DQN.     Convex Q-learning was the winner in terms of reliability and performance for reasons that are not clear at this time.  
	\notes{
		We wrote `` in the numerical examples reported here'' but provide no evidence.  A comparison is needed for thesis
		\\
		\fl{ok.}	
}

Theory concerning the impact of disturbances and measurement noise  is an important area for future research.  We believe the value of filtering in convex Q-learning will be apparent when we include measurement noise in simulation experiments, but currently have no guidelines to optimize $\sigma$ or opt for a different approach to smoothing the Q-ODE.
	
	The use of state space collapse to design a function approximation architecture was very successful in the example considered.   This will likely prove valuable in other   applications.    Such extensions may require techniques to characterize or approximate the manifold on which an optimal solution evolves, or perhaps we can create algorithms that will ``learn'' this structure.

	


	\appendix
	
	\section{Appendix}
	The proof of \Cref{t:convexQbdd} requires  Grönwall's inequality in this simplified form:
	\begin{lemma}[Bellman-Grönwall]   
		\label[lemma]{t:B-GI}
		Let  $w$ be a continuous real-valued function on the interval $[0,\Hor]$. 
		Suppose that the following integral bound holds with the constants $\alpha , \beta \ge 0$:
	\[
		w_r\leq \alpha +   \beta \int_0^r  w_s\, ds    \,,\qquad 0\le r\le T
	\]
		Then, $ w_r\leq \alpha   e^{\beta r} $ for $ 0\le r\le T$.
	\end{lemma}
	
	\ProofOf{Proof of \Cref{t:convexQbdd}}
	Since $H \ge \uH$, it follows from \eqref{e:q_integral_lower} that for any trajectory, 
	\begin{align}
		\begin{split}
			\Hcev_r  & \ge -\sigma e^{ -\sigma r}   J_0(x_\Hor) + \clHcev_r - \clCcev_r.
		\end{split}
		\label{e:q_integral_lower_mod}
	\end{align} 
	For any optimal trajectory $\{x^\star_r, u^\star_r\}$  
	we have $\Hcev_r^\star  = \uHcev_r^\star$, so from \eqref{e:q_algebraic}, 
	\begin{align}
		\begin{split}
			\Hcev^\star &= -\sigma e^{ -\sigma r}   J_0(x_\Hor) 
			+ \clHcev^\star_r - \clCcev^\star_r
		\end{split}
		\label{e:q_integral_mod}
	\end{align}
	
Denote 
$
	\Delta_r \eqdef   \Hcev_r^\star - \Hcev_r.
$
Subtracting   \eqref{e:q_integral_lower_mod} from \eqref{e:q_integral_mod} then yields,
	\begin{equation*}
		\Delta_r \le \clHcev^\star_r - \clHcev_r = \int_0^r  e^{-\sigma (r-s)} \Delta_s \, ds,
	\end{equation*}
	where the equality on the right follows from the definitions of $\clH$ and $\clH^\star$.
	Setting $w_r = e^{\sigma r} \Delta_r$ and applying \Cref{t:B-GI} gives
\[
	w_r = e^{\sigma r} \Delta_r \le 0, \qquad 0 \le r \le \Hor,
\]
	which in turn implies $\Delta_r \le 0$, thereby yielding $H \ge H^\star$ along this optimal trajectory.
	
	It follows that $H(x,u,r) \ge 
	H^\star(x,u,r)$ for any $(x,u, r)$, since there is an optimizing trajectory that passes through any such triple. 
	\qed
	
	\ProofOf{Proof of \Cref{t:AssRel}}
	To establish that   Condition E1 implies E2, we establish the contrapositive: if there is a non-zero vector $v$ satisfying $\uHcev^v_r \ge \clHcev^v_r$ for each $r$, then the set $\{ \psicev_r: 0 \le r \le \Hor \}$ is restricted to a half space in $\Re^d$.
	
If such $v$ exists, then by definition of $\uH$,  
\[
	H^v(x_{\Hor - r}, u, r) \ge \uHcev^v_r \ge \clHcev^v_r \,, \qquad u \in \Re^m.
\]
Letting $p_r = v^\transpose \psicev_r$, and $y_r = v^\transpose \Psicev_r$,   this inequality implies that $p_r \ge y_r$ and by definition, 
	\begin{equation}
		\begin{aligned}  			
			y_r &  = \sigma \int_{0}^r e^{-\sigma(t-r)} p_\tau d\tau
			\\
			\frac{d}{dr} y_r & = -\sigma (y_r - p_r), 
		\end{aligned} 
		\label{eq:intFom}	
	\end{equation}
On applying the  boundary condition $y_0 = 0$,  
\[
	y_r = -\sigma \int_{0}^r (y_\tau - p_\tau) d\tau \ge 0.
\]
	
	Letting $ \delta_r =	p_r - y_r$,  which is non-negative,    gives
	$p_r =  y_r + \delta_r$,  and  for each  $0\le r \le \Hor$,
	\begin{equation*}
		\begin{aligned}
			v^\transpose \psicev_r  = -\sigma \int_{0}^r (y_\tau - p_\tau) d\tau + \delta_r 
			\ge 0  \,, .
		\end{aligned}
	\end{equation*}
	Hence Condition E1 fails when   E2 fails, as claimed.   
	\notes{
		sm2fl:  you missed one thing.   $\{\psi_t\}$ is restricted to a half space, is the same thing as  $\{\psi_{\Hor-r}\}$ is restricted to a half space.
	\\
	\fl{Got it.}
}
	
	To show that   Condition E1 implies E3, we again establish the contrapositive: if $\det(\Sigma) = 0$, then the set $\{ \psi_t: 0 \le t\le \Hor \}$ is restricted to a half space in $\Re^d$.
	
If $v \in \text{Null}( \Sigma )$ with $v\not=0$, then 
\[
	0 = v^\transpose \Sigma v = \frac{1}{\Hor} \int_0^\Hor (v^\transpose \tilpsi_t)^2  \, dt
\] 
	Since $\{ \tilpsi_t \}$ is continuous in $t$, it follows that
\[
	v^\transpose \tilpsi _t = 0 \,, \qquad \text{for $0 \le t \le \Hor$}.
\]
	This implies that $\{\psi_t\}$ is restricted to a half space, so that Condition E1 fails.
	\qed

\ProofOf{Proof of \Cref{t:AssConclusions}}
There are two parts to the proof.  We first establish that  $\Uptheta$ is bounded under E1.    \Cref{t:AssRel}  tells us that E2 follows from E1,   so  it suffices to show that if Condition E2 holds then $\Uptheta$ is bounded.   We establish its contrapositive: if $\Uptheta$ is unbounded, then there is a non-zero vector $v$ satisfying $\uHcev^v_r \ge \clHcev^v_r$ for $0 \le r \le \Hor$.
	
	If $\Uptheta$ is unbounded, then for each $m\ge 0$, there exists $\theta^{m}$ such that $\| \theta^{m} \| \ge m$, and
	\begin{align}
		\label{e:intolp}
		\frac{1}{\Hor} \int_{0}^{\Hor} \max\bigg\{0, \clJ_r-\uHcev^{\theta^m}_r + \clHcev^{\theta^m}_r\bigg\}\, dr \le \Tol
	\end{align}
	with $\clJ_r \eqdef -e^{-\sigma r}J_0(x_\Hor) - \clCcev_r$.
	
	Dividing \eqref{e:intolp} by $\|\theta^{m}\|$ gives:
	\begin{equation}
		\label{e:infineq}
		\frac{1}{\Hor} \int_{0}^{\Hor} \max\bigg\{0,
		\tfrac{\clJ_r}{\|\theta^{m}\|}
		-\tfrac{\uHcev^{\theta^{m}}_r}{\|\theta^{m}\|}
		+ \tfrac{\clHcev^{\theta^{m}}_r}{\|\theta^{m}\|}
		\bigg\} dr
		\le \frac{\Tol}{\|\theta^{m}\|}
	\end{equation}
	Denote $\ctheta^m =  {\theta^m}/{\| \theta^m\|}$. By the definition of $\uHcev^{\theta^m}$,  
	\begin{align*}
		\frac{1}{\| \theta^m \|} \uHcev^{\theta^m}_r &=  \min_u  \Bigl\{ \frac{1}{\| \theta^m \|} H^{\theta^m}(\xcev_r, u, r)  \Bigr\} = \uHcev^{\ctheta^m}_r
	\end{align*}
	Thus, we can write \eqref{e:infineq}  as
	\begin{equation}
		\label{e:infeqde}
		\frac{1}{\Hor} \int_{0}^{\Hor} \max\bigg\{0, \frac{\clJ_r}{\|\theta^m\|}
		-\uHcev^{\ctheta^{m}}_r 
		+ \clHcev^{\ctheta^{m}}_r
		\bigg\}dr
		\le \frac{\Tol}{\|\theta^{m}\|}
	\end{equation}
	Since $\|\ctheta^m\|  = 1$ for each $m$, there exists a convergent subsequence $\{ \theta^{m_i} \}$ with limit satisfying  $\| \ctheta  \| = 1$:
\[
	\lim_{i\to\infty} \frac{ \theta^{m_i} }{ \| \theta^{m_i} \| } = \lim_{i\to\infty} \ctheta^{m_i} = \ctheta
\]
	The inequality  \eqref{e:infeqde} then gives
	\begin{equation*}
		\begin{aligned}
			\frac{1}{\Hor} &\int_{0}^{\Hor}  \max\Big\{0, -\uHcev^{\ctheta}_r  
			+ \clHcev^{\ctheta}_r
			\Big\}\, dr
			\\
			&=\lim_{i\to\infty}\frac{1}{\Hor} \int_{0}^{\Hor} \max\Big\{0, 
			\tfrac{1}{\|\theta^{m_i}\|}\clJ_r
			-\uHcev^{\ctheta^{m_i}}_r 
			+ \clHcev^{\ctheta^{m_i}}_r
			\Big\} \, dr
			\\
			&\le 0
		\end{aligned}
\end{equation*}
Continuity of $\{\uHcev^{\ctheta}_r  , \clHcev^{\ctheta}_r : 0\le r\le \Hor\}$	implies   the desired conclusion: E2 fails, with $v= \ctheta$,
\[
	\uHcev^{\ctheta}_r \ge \clHcev^{\ctheta}_r \,, \qquad 0\le r \le \Hor \,.
\]

For the converse we once again establish the contrapositive:   if Condition E2 fails, we show that $\Uptheta$ is unbounded.  

Failure of E2 implies that there is   $v\neq 0$ satisfying $\uHcev^v_r \ge \Hcev^v_r$ for $0 \le r \le \Hor$.   To show that $\Uptheta$ is  unbounded we fix $\theta^0 \in \Uptheta$,  and show that $ \theta^\omega \eqdef \theta^0 + \omega v  \in \Uptheta$ for each $\omega\ge 0$.
	Because the function class is linear, we have
	\begin{align*}
		\uHcev^{\theta^{\omega}}_r 
		&\eqdef \min_u \{ H^{\theta^{\omega}}(x_{\Hor - r}, u,  \Hor - r)  \}
		\\
		&= \min_u \{ H^{\theta^{0}}(x_{\Hor - r}, u,  \Hor - r) + \omega   H^{v}(x_{\Hor - r}, u,  \Hor - r)\}
	\end{align*}
	This and sub-linearity of the minimum gives for each $r$,
	\begin{align*}
		\uHcev^{\theta^{\omega}}_r  \ge     \uHcev^{\theta^{0}}_r  +      \omega\uHcev^{v}_r.
	\end{align*}
	It follows that the Bellman error for   $\theta^{\omega}$ admits the bound,
	\begin{align*}
		\clB^{\theta^\omega}_r &\eqdef
		\clJ_r-\uHcev^{\theta^{\omega}}_r + \clHcev^{\theta^{\omega}}_r
		\\
		&\le \clJ_r-[\uHcev^{\theta^{0}}_r + \omega\uHcev^{v}_r] + [\clHcev^{\theta^{0}}_r + \omega\clHcev^{v}_r]
	\end{align*}
	and on rearranging terms,  
\[
	\clB^{\theta^\omega}_r 
	\le \clB^{\theta^0}_r + \omega[-\uHcev^{v}_r + \clHcev^{v}_r]
\]
	
	By  assumption, we have $-\uHcev^{v}_r + \clHcev^{v}_r \le 0$ and thus
	$\clB^{\theta^\omega}_r \le \clB^{\theta^0}_r$.   Consequently,
		   $\theta^\omega \in\Uptheta$ for every $\omega$, as claimed:
$$
	\begin{aligned}
		\frac{1}{\Hor} \int_0^\Hor \max\{ 0, \clB^{\theta^\omega}_r \}\, dr 
		&\le \frac{1}{\Hor} \int_0^\Hor \max\{ 0, \clB^{\theta^0}_r \}\, dr \le \Tol  
	\end{aligned}
	\qed
$$

	
	\bibliographystyle{abbrv}
	\bibliography{strings,markov,q,CollapseExtras,PolicyCollapseExtras}  

	{}
	\null  
	\null  

\end{document}